\documentclass[12pt]{article}
\usepackage{booktabs}
\usepackage{caption}
\usepackage{mathrsfs}
\usepackage{amsmath}
\usepackage{amsfonts,amsthm,amssymb,mathrsfs,bbding}
\usepackage{float}
\usepackage{graphicx}
\usepackage{enumerate}
\usepackage{tikz}
\usepackage{caption}
\usepackage{rotating}
\allowdisplaybreaks[4]
\usepackage{tabularx}
\usepackage{cite}
\evensidemargin=\oddsidemargin\topmargin=-1.5cm

\newtheorem{thm}{Theorem}[section]
\newtheorem{prob}{Problem}

\newtheorem{lem}[thm]{Lemma}
\newtheorem{cor}[thm]{Corollary}

\theoremstyle{definition}
\renewcommand\proofname{\bf Proof}
\addtocounter{section}{0}
\begin{document}

\title{\bf Spectral conditions for $k$-extendability and $k$-factors of  bipartite graphs\footnote{This work is supported by the National Natural Science Foundation of China (Grant Nos. 12271162 and 12011530064) and Natural Science Foundation of Shanghai (No. 22ZR1416300).}}
\author{Dandan Fan$^{a,b}$, Huiqiu Lin$^a$\thanks{Corresponding author.~~Email address: huiqiulin@126.com (H. Lin).}\\[2mm]
\small\it $^a$Department of Mathematics, East China University of Science and Technology, \\
\small\it   Shanghai 200237, P.R. China\\[1mm]
\small\it $^b$College of Mathematics and Physics, Xinjiang Agricultural University,\\
\small\it Urumqi, Xinjiang 830052, China}
\date{}
\maketitle
{\flushleft\large\bf Abstract}
 Let $G$ be a connected graph. If $G$ contains a matching of size $k$, and every matching of size $k$ is contained in a perfect matching of $G$, then $G$ is said to be \emph{$k$-extendable}. A $k$-regular spanning subgraph of $G$ is called a   \textit{$k$-factor}. In this paper, we provide spectral conditions for a (balanced bipartite) graph with minimum degree $\delta$ to be $k$-extendable, and for the existence of a $k$-factor in a balanced bipartite graph, respectively. Our results generalize some previous results on perfect matchings of graphs, and extend the results in \cite{D.F} and \cite{W.L} to $k$-extendable graphs. Furthermore, our results generalize the result of Lu, Liu and Tian \cite{Lu-Liu} to general regular factors. Additionally, using the equivalence of $k$ edge-disjoint perfect matchings and $k$-factors in balanced bipartite graphs, our results can derive a spectral condition for the existence of $k$ edge-disjoint perfect matchings in balanced bipartite graphs.

\begin{flushleft}
\end{flushleft}

\section{Introduction}
Perfect matchings theory, which studies the simplest nontrivial substructures of graphs, is one of the earliest reasearch areas in graph theory. The characterization of perfect matchings was initially given by Frobenius\cite{Frobenius} in 1917, who proved that a bipartite graph of order $n$ has a perfect matching if and only if each vertex cover has size at least $n/2$. Soon afterwards, Hall\cite{Hall} and Tutte\cite{Tutte} provided better and simpler criterions for the existence of perfect matchings in
bipartite graphs and general graphs, respectively. From then on, many researchers have attempted to find  sufficient conditions for the existences of perfect matchings by using various graph parameters (see \cite{Petersen,Sumner-1,Win}). In order to enrich the results of the graph structures with perfect matchings, one stream concerned perfect matching which has attracted a lot of interest is that of researching the structures and properties of graphs containing many perfect matchings, such as $k$-factor-critical graphs (see \cite{Favaron} for the definition), $k$-extendable graphs and other various matching-extendable graphs. In 1979, Sumner\cite{Sumner} asked to characterize the graphs with the property that every matching can be extended to a perfect matching.  It turns out that only complete graphs of even order and complete balanced bipartite graphs satisfy this property. By relaxing the property a bit, Plummer\cite{Plummer} introduced the concept of $k$-extendable graph by requiring only the matching with the same size to be extended to a perfect matching. Up to now, much attention has been paid on  various graph parameters of $k$-extendable graphs, such as independence number \cite{Maschlanka}, connectivity\cite{Plummer-3}, toughness\cite{Plummer-2} and genus\cite{Plummer-1}. For more comprehensive details, readers can refer to \cite{L-P}.

% Another problem concerned perfect matching is extending it to a general factor.

One of the main goals of this paper is to investigate the $k$-extendability of graphs from a spectral perspective. It is easy to see that if $G$ is $0$-extendable, then $G$ contains a perfect matching. In the past few years, many researchers focused on finding the relationship between the eigenvalues and perfect matchings of a graph. In 2005, Brouwer and Haemers \cite{A.B}  initially  described a regular graph to contain a perfect matching in terms of the third largest adjacency eigenvalue, which was improved in \cite{S.C-1,S.C,S.C-2}. Recently, O\cite{S.O} gave a spectral condition to guarantee the existence of a perfect matching in a connected graph. By imposing the minimum degree of a graph as a parameter, Liu, Liu and Feng \cite{W.L} and Fan et al.\cite{D.F} extended the result of O \cite{S.O} in graphs and bipartite graphs, respectively. In this paper, we generalize their results by giving a spectral condition to guarantee a (balanced bipartite) graph with minimum degree $\delta$ to be $k$-extendable.

Let $G$ be a graph with adjacency matrix $A(G)$. The largest eigenvalue of $A(G)$, denoted by $\rho(G)$, is called the \textit{spectral radius} of $G$. Denote by `$\nabla$' and `$\cup$' the join and union products, respectively.
Let $F(k,\delta)=\max\{8\delta-10k+4, \delta(\delta-2k)^2+\delta-1\}$.
\begin{thm}\label{thm::1.1}
Suppose that $G$ is a connected graph of even order $n\geq F(k,\delta)$ with minimum degree $\delta\geq 2k$, where $k\geq 1$. If $\rho(G)\geq\rho(K_{\delta}\nabla(K_{n-2\delta+2k-1}\cup(\delta-2k+1)K_1))$, then $G$ is $k$-extendable, unless $G\cong K_{\delta}\nabla(K_{n-2\delta+2k-1}\cup(\delta-2k+1)K_1)$.
\end{thm}

Given two bipartite graphs $G_1=(X_{1}, Y_{1})$ and $G_2=(X_{2}, Y_{2})$, let $G_{1}\nabla_{1} G_2$ denote the graph obtained from $G_1 \cup G_2$ by adding all possible edges between $X_{1}$ and $Y_{2}$. A bipartite graph $G=(X,Y)$ is called \emph{balanced} if $|X|=|Y|$. Clearly,  every $k$-extendable bipartite graph is balanced.
\begin{thm}\label{thm::1.2}
Let $1\leq k\leq \frac{n}{2}-1$, and let $G$ be a balanced bipartite graph of order $n$ with minimum degree $\delta$. If
$$\rho(G)\geq \rho(K_{\delta,\delta+k+1}\nabla_{1}K_{\frac{n}{2}-\delta,\frac{n}{2}-\delta-k-1}),$$
then $G$ is $k$-extendable, unless $G\cong K_{\delta,\delta+k+1}\nabla_{1}K_{\frac{n}{2}-\delta,\frac{n}{2}-\delta-k-1}$.
\end{thm}

A $k$-regular spanning subgraph of a graph  is called a \textit{$k$-factor}. In particular, a perfect matching is a $1$-factor. Based on the classic works of Tutte\cite{Tutte} and Hall\cite{L-P}, many researchers have made efforts to find structural conditions for the existence of $k$-factors in graphs  \cite{Rado,Enomoto,Katerinis,Ore-1,Lu-Ning}. The existence of $k$-factors has also been well studied from spectral perspectives.
A hamiltonian cycle is a connected $2$-factor of a graph. In 2010, Fiedler and Nikiforov\cite{Fiedler} gave a spectral condition for the existence of a hamiltonian cycle in a graph, and this result was improved in \cite{Ge-Ning,Ge-Ning-1, Chen-1}. Extending the above results to general factors, Cho, Hyun, O and Park\cite{E.C} conjectured a spectral condition for the existence of $k$-factors in graphs. Subsequently, Fan, Lin and Lu\cite{F.D-1} confirmed this conjecture for $n\geq 3k-1$. For balanced bipartite graphs, Lu, Liu and Tian\cite{Lu-Liu} investigated the spectral conditions for the existence of a hamiltonian cycle in complementary graphs. Motivated by the recent works in \cite{F.D-1} and \cite{Lu-Liu}, we give a spectral condition for the existence of a $k$-factor in a balanced bipartite graph. Denote by $G\backslash E(H)$ the graph obtained from $G$ by deleting the edges of $H$, where $H$ is a subgraph of $G$.
\begin{thm}\label{thm::1.3}
Let $2\leq k\leq \frac{n}{2}-1$, and let $G$ be a connected balanced bipartite graph with order $n$. If
$$\rho(G)\geq \rho(K_{\frac{n}{2},\frac{n}{2}}\backslash E(K_{1,\frac{n}{2}-k+1})),$$
then $G$ contains a $k$-factor, unless $G\cong K_{\frac{n}{2},\frac{n}{2}}\backslash E(K_{1,\frac{n}{2}-k+1})$.
\end{thm}

In 2008, Liu et al.\cite{L-G} proved the equivalence of $k$ edge-disjoint perfect matchings and $k$-factors in balanced bipartite graphs. By Theorem \ref{thm::1.3}, we immediately deduce the following result.

\begin{cor}
Let $2\leq k\leq \frac{n}{2}-1$, and let $G$ be a connected balanced bipartite graph with order $n$. If
$$\rho(G)\geq \rho(K_{\frac{n}{2},\frac{n}{2}}\backslash E(K_{1,\frac{n}{2}-k+1})),$$
then $G$ contains $k$ edge-disjoint perfect matchings, unless $G\cong K_{\frac{n}{2},\frac{n}{2}}\backslash E(K_{1,\frac{n}{2}-k+1})$.
\end{cor}
%\begin{remark}
%A hamiltonian cycle is a connected 2-factor of a graph.
%\end{remark}

\section{Proof of Theorems~1.1 and 1.2}
\indent In this section, we present the proofs of Theorems~\ref{thm::1.1} and \ref{thm::1.2}. In \cite{Chen}, Chen established a sufficient and necessary condition for a graph to be $k$-extendable. For any $S\subseteq V(G)$, let $G[S]$ be the subgraph of $G$ induced by $S$.
\begin{lem}[Chen\cite{Chen}]\label{lem::2.1}
For $k\geq 1$, a graph $G$ is $k$-extendable if and only if
$$o(G-S)\leq |S|-2k,$$
for any $S\subseteq V(G)$ such that $G[S]$ contains $k$ independent edges, where $o(H)$ is the number of odd components in a graph $H$.
\end{lem}

\begin{lem}[See \cite{D.F}]\label{lem::2.2}
Let $n=\sum_{i=1}^t n_i+s$. If $n_{1}\geq n_{2}\geq \cdots\geq n_{t}\geq p$ and $n_{1}<n-s-p(t-1)$, then
$$\rho(K_{s} \nabla (K_{n_{1}}\cup K_{n_{2}}\cup \cdots \cup K_{n_{t}}))<\rho(K_{s} \nabla (K_{n-s-p(t-1)}\cup (t-1)K_{p})).$$
\end{lem}

\begin{lem}\label{lem::2.3}
Let $n\geq 8\delta-10k+4$, where $k\geq 1$ and $\delta\geq 2k+1$. Then
$$\rho(K_{2k}\nabla(K_{\delta-2k+1}\cup K_{n-\delta-1}))<\rho(K_{\delta}\nabla(K_{n-2\delta+2k-1}\cup (\delta-2k+1)K_{1})).$$
\end{lem}
\renewcommand\proofname{\bf Proof}
\begin{proof}
Let $G=K_{2k}\nabla(K_{\delta-2k+1}\cup K_{n-\delta-1})$. Then the vertex set of $G$ can be partitioned as $V(G)=V(K_{2k})\cup V(K_{\delta-2k+1})\cup V(K_{n-\delta-1})$, where $V(K_{2k})=\{v_1,\ldots,v_{2k}\}$, $V(K_{\delta-2k+1})=\{u_{1},\ldots,u_{\delta-2k+1}\}$ and $V(K_{n-\delta-1})=\{w_{1},\ldots,w_{n-\delta-1}\}$. Suppose that $E_{1}=\{u_{i}w_{j}|~1\leq i\leq \delta-2k+1, 1\leq j\leq \delta-2k\}$ and $E_{2}=\{u_{i}u_{j}|~
 1\leq i\leq \delta-2k, i+1\leq j\leq \delta\!-2k+\!1\}$.  Let $G'=G+E_1-E_2$. Clearly, $G'\cong K_{\delta}\nabla(K_{n-2\delta+2k-1}\cup (\delta-2k+1)K_{1})$.
  Let $x$ be the Perron vector of $A(G)$. By symmetry, $x$ takes the same value (say $x_1$, $x_2$ and $x_3$) on the vertices of $V(K_{2k})$, $V(K_{\delta-2k+1})$ and $V(K_{n-\delta-1})$, respectively. Then, by $A(G)x=\rho(G) x$, we have
$$
\left\{
\begin{aligned}
 &\rho x_2=2k x_{1}+(\delta-2k)x_2,\\
 &\rho x_3=2k x_{1}+(n-\delta-2) x_{3},
  \end{aligned}
\right.
$$
which leads to
$$(\rho-\delta+2k)(x_3-x_2)=(n-2\delta+2k-2)x_3>0$$
because $n\geq 8\delta\!-\!10k\!+\!4$, $k\geq 1$ and $\delta\geq 2k$. Note that
$K_{n-\delta+2k-1}\cup K_{\delta-2k+1}$ is a proper subgraph of $G$. Then $\rho(G)>n-\delta+2k-2>\delta-2k$. Therefore,
$$x_3>x_2,$$
and hence
\begin{equation*}
\begin{aligned}
x^{T}(A(G')-A(G))x&\geq 2\left(\sum_{i=1}^{\delta-2k}\sum_{j=1}^{\delta-2k+1}x_{w_i}x_{u_j}-
\sum_{i=1}^{\delta-2k}\sum_{j=i+1}^{\delta-2k+1}x_{u_i}x_{u_j}\right)\\
&=2x_2(\delta-2k)(\delta-2k+1)\left(x_{3}-\frac{x_2}{2}\right)\\
&>0.
\end{aligned}
\end{equation*}
The result follows.
\end{proof}

\begin{lem}[Plummer \cite{Plummer}]\label{lem::2.4}
A bipartite graph $G=(A,B)$ is $k$-extendable if and only if $|A|=|B|$ and
$$|N(X)|\geq |X|+k$$
for all nonempty subset $X\subseteq A$ with $|X|\leq |A|-k$.
\end{lem}

Let $M$ be a real $n\times n$ matrix, and let $X=\{1,2,\ldots,n\}$. Given a partition $\Pi:X=X_{1}\cup X_{2}\cup \cdots \cup X_{k}$,  the matrix $M$ can be correspondingly partitioned as
$$
M=\left(\begin{array}{ccccccc}
M_{1,1}&M_{1,2}&\cdots &M_{1,k}\\
M_{2,1}&M_{2,2}&\cdots &M_{2,k}\\
\vdots& \vdots& \ddots& \vdots\\
M_{k,1}&M_{k,2}&\cdots &M_{k,k}\\
\end{array}\right).
$$
The \textit{quotient matrix} of $M$ with respect to $\Pi$ is defined as the $k\times k$ matrix $B_\Pi=(b_{i,j})_{i,j=1}^k$ where $b_{i,j}$ is the  average value of all row sums of $M_{i,j}$.
The partition $\Pi$ is called \textit{equitable} if each block $M_{i,j}$ of $M$ has constant row sum $b_{i,j}$.
Also, we say that the quotient matrix $B_\Pi$ is \textit{equitable} if $\Pi$ is an equitable partition of $M$.

\begin{lem}[Brouwer and Haemers \cite{BH}, p. 30; Godsil and Royle \cite{C.Godsil}, pp. 196--198]\label{lem::2.5}
Let $M$ be a real symmetric matrix, and let $\lambda_{1}(M)$ be the largest eigenvalue of $M$. If $B_\Pi$ is an equitable quotient matrix of $M$, then the eigenvalues of  $B_\Pi$ are also eigenvalues of $M$. Furthermore, if $M$ is nonnegative and irreducible, then $\lambda_{1}(M) = \lambda_{1}(B_\Pi).$
\end{lem}

\begin{lem} \label{lem::2.6}
For $n\geq 4s+2k+2$ and $k\geq 1$, we have
$$\rho(K_{s,s+k+1}\nabla_{1} K_{\frac{n}{2}-s,\frac{n}{2}-s-k-1})<\rho(K_{s-1,s+k}\nabla_{1} K_{\frac{n}{2}-s+1,\frac{n}{2}-s-k}).$$
\end{lem}

\renewcommand\proofname{\bf Proof}
\begin{proof}
Note that $A(K_{s,s+k+1}\nabla_{1} K_{\frac{n}{2}-s,\frac{n}{2}-s-k-1})$ has the equitable quotient matrix
$$
B_\Pi^s=\begin{bmatrix}
0 &0 & s+k+1& \frac{n}{2}-s-k-1\\
0&0&0 &\frac{n}{2}-s-k-1\\
  s&0& 0& 0\\
   s &\frac{n}{2}-s &0 &0
\end{bmatrix}.
$$
By a simple computation, the characteristic polynomial of $B_\Pi^s$ is
\begin{equation*}
\begin{aligned}
\varphi(B_{\Pi}^s,x)\!=\!x^4\!+\!\left(\frac{(2k\!+\!2s\!+\!2\!-\!n)n}{4}\!-\!(s\!+\!k\!+\!1)s\right)x^2\!-\!\frac{s(n\!-\!2s)(s\!+\!k\!+\!1)(2s\!-\!n\!+\!2k\!+\!2)}{4}
\end{aligned}
\end{equation*}
Observe that $A(K_{s-1,s+k}\nabla_{1} K_{\frac{n}{2}-s+1,\frac{n}{2}-s-k})$ has the equitable quotient matrix $B_{\Pi}^{s-1}$, which is obtained by replacing $s$ with $s-1$ in $B_\Pi^s$. Then
\begin{equation*}
\begin{aligned}
\varphi(B_{\Pi}^{s},x)-\varphi(B_{\Pi}^{s-1},x)&=\frac{1}{4}(n\!-\!2k\!-\!4s)((k\!+\!2s)n\!-\!4ks\!+\!2x^2\!-\!4s^2\!+\!2k)\\
&\geq \frac{1}{2}(2k^2\!+\!(4s\!+\!4)k\!+\!4s^2\!+\!2x^2\!+\!4s)~~(\mbox{since $n\geq 4s\!+\!2k\!+\!2$})\\
&>0,
\end{aligned}
\end{equation*}
which leads to
$\lambda_{1}(B_{\Pi}^{s})<\lambda_{1}(B_{\Pi}^{s-1})$. It follows that
$$\rho(K_{s,s+k+1}\nabla_{1} K_{\frac{n}{2}-s,\frac{n}{2}-s-k-1})<\rho(K_{s-1,s+k}\nabla_{1} K_{\frac{n}{2}-s+1,\frac{n}{2}-s-k})$$
by Lemma \ref{lem::2.5}.

This completes the proof.
\end{proof}

Now we give the proof of Theorem \ref{thm::1.1}.
\renewcommand\proofname{\bf Proof of Theorem \ref{thm::1.1}}
\begin{proof}

Suppose that $G$ is not $k$-extendable, by Lemma \ref{lem::2.1}, there exists some nonempty subset $S$ of $V(G)$ such that $s=|S|\geq 2k$ and $o(G-S)>s-2k$.
Since $n$ is even, $o(G-S)$ and $s$ have the same parity,
we have $o(G-S)\geq s-2k+2$. It is clear that $G$ is a spanning subgraph of $G^1=K_{s} \nabla (K_{n_1}\cup K_{n_2}\cup \cdots \cup K_{n_{s-2k+2}})$ for some positive odd integers $n_1\geq n_2\geq \cdots\geq n_{s-2k+2}$ with $\sum_{i=1}^{s-2k+2}n_i=n-s$. Thus,
\begin{equation}\label{equ::1}
\rho(G)\leq\rho(G^1),
\end{equation}
where the equality holds if and only if $G\cong G^1$. We shall divide the proof into the following three cases.

{\flushleft\bf{Case 1.}} $s\geq\delta+1.$
\\[2mm]
\indent Let $G^2=K_{s} \nabla (K_{n-2s+2k-1}\cup (s-2k+1)K_{1})$. Note that $s\geq \delta+1$. Then Lemma \ref{lem::2.2} gives that
\begin{equation}\label{equ::2}
\rho(G^1)\leq \rho(G^2),
\end{equation}
where the equality holds if and only if $(n_1,\ldots,n_{s-2k+2})=(n-2s+2k-1,1,\ldots,1)$.
The vertex set of $G^2$ can be partitioned as $V(G^2)=V(K_{s})\cup V((s-2k+1)K_{1})\cup V(K_{n-2s+2k-1})$, where $V(K_s)=\{v_1,\ldots,v_s\}$, $V((s-2k+1)K_{1})=\{u_{1},\ldots,u_{s-2k+1}\}$ and $V(K_{n-2s+2k-1})=\{w_{1},\ldots,w_{n-2s+2k-1}\}$. Suppose that $E_{1}=\{u_{i}w_{j}|~\delta\!-\!2k+2\leq i\leq s\!-\!2k\!+\!1, 1\leq j\leq n-2s+2k-1\}\cup\{u_{i}u_{j}|~ \delta\!-\!2k+2\leq i\leq s-2k, i\!+\!1\leq j\leq s-2k\!+\!1\}\!$ and $E_{2}=\{v_{i}u_{j}|~\delta\!+\!
 1\leq i\leq s, 1\leq j\leq \delta\!-2k+\!1\}$.  Let $G' =G^2\!+E_{1}-E_{2}$.
Clearly, $G'\cong K_{\delta} \nabla (K_{n-2\delta+2k-1}\cup (\delta-2k+1)K_{1})$.
Let $x$ be the Perron vector of $A(G^2)$, and let $\rho=\rho(G^2)$. By symmetry, $x$ takes the same value (say $x_1$, $x_2$ and $x_3$) on the vertices of $V(K_s)$, $V((s-2k+1)K_{1})$ and $V(K_{n-2s+2k-1})$, respectively. Then, by $A(G^2)x=\rho x$, we have
\begin{eqnarray*}
 &\rho x_2&=s x_{1},\\
 &\rho x_3 &=s x_{1}+(n-2s+2k-2) x_{3}.
\end{eqnarray*}
Observe that $n\geq 2s-2k+2$. Then $x_{3}\geq x_{2}$ and
\begin{eqnarray}
 x_{2}=\frac{sx_1}{\rho}. \label{equ::3}
\end{eqnarray}
Similarly, let $y$ be the Perron vector of $A(G')$, and let $\rho(G')=\rho'$. By symmetry, $y$ takes the same values $y_{1}$, $y_{2}$ and $y_{3}$ on the vertices of $V(K_\delta)$, $V((\delta-2k+1)K_{1})$ and $V(K_{n-2\delta+2k-1})$, respectively. Then, by $A(G')y=\rho' y$, we have
\begin{eqnarray}
 &\rho' y_2&=\delta y_{1},\label{equ::4}\\
 &\rho' y_3 &=\delta y_{1}+(n-2\delta+2k-2) y_{3}.\label{equ::5}
\end{eqnarray}
Note that $G'$ contains $K_{n-2\delta+2k-1}\cup K_{\delta}\cup (\delta-2k+1)K_1$ as a proper  subgraph. Then $\rho'>n-2\delta+2k-2$. Combining this with (\ref{equ::4}) and (\ref{equ::5}), we have
\begin{equation} \label{equ::6}
      y_3=\frac{\rho' y_2}{\rho'-(n-2\delta+2k-2)}.
\end{equation}
Note that $n\geq 2s-2k+2$. Then $\delta+1\leq s\leq(n+2k-2)/2$. As $G^2$ is not a complete graph, $\rho<n-1$. Now we shall prove that  $\rho<\rho'$.  Indeed, if  $\rho\geq\rho'$, then from  (\ref{equ::3}) and (\ref{equ::6}), we obtain
\begin{equation*}
\begin{aligned}
   & y^{T}(\rho'-\rho)x\\
   &=y^{T}(A(G')-A(G^2))x\\
   &=\sum_{i=\delta\!-\!2k\!+\!2}^{s\!-\!2k\!+\!1}\sum_{j=1}^{n\!-\!2s\!+\!2k\!-\!1}\!(x_{u_{i}}y_{w_{j}}\!+\!x_{w_{j}}y_{u_{i}})\!
+\sum_{i=\delta\!-\!2k\!+\!2}^{s-2k}\sum_{j=i\!+\!1}^{s\!-\!2k\!+\!1}\!(x_{u_{i}}y_{u_{j}}\!+\!x_{u_{j}}y_{u_{i}})\! \!-\!
\sum_{i=\delta\!+\!1}^{s}\sum_{j=1}^{\delta\!-\!2k\!+\!1}\!(x_{v_{i}}y_{u_{j}}\!+\!x_{u_{j}}y_{v_{i}})\!\\
   &=(s\!-\!\delta)[(n\!-\!2s\!+\!2k\!-\!1)(x_{2}y_{3}\!+\!x_{3}y_{3})\!+\!(s\!-\!\delta\!-\!1)x_{2}y_{3}\!-\!(\delta\!-\!2k\!+\!1)(x_{1}y_{2}\!+\!x_{2}y_{3})]\\
   &=(s\!-\!\delta)[(n\!-\!2\delta\!-\!s\!+\!4k\!-\!3)x_{2}y_{3}\!+\!(n\!-\!2s\!+\!2k\!-\!1)x_{3}y_{3}\!-\!(\delta\!-\!2k\!+\!1)x_{1}y_{2}]\\
   &\geq(s\!-\!\delta)[(2n\!-\!2\delta\!-\!3s\!+\!6k\!-\!4)x_{2}y_{3}\!-\!(\delta\!-\!2k\!+\!1)x_{1}y_{2}]  ~(\mbox{since $x_{3}\!\geq\! x_{2}$})\\
   &=(s\!-\!\delta)x_{1}y_{2}\left(\frac{s\rho'(2n\!-\!2\delta\!-\!3s\!+\!6k\!-\!4)}{\rho(\rho'-(n-2\delta+2k-2))}-(\delta-2k+1)\right)
   ~~(\mbox{by (\ref{equ::3}) and (\ref{equ::6})})\\
   &>\frac{(\delta\!-\!2k\!+\!1)(s\!-\!\delta)x_{1}y_{2}}{\rho(\rho'\!-\!(n\!-\!2\delta\!+\!2k\!-\!2))}[\rho'(2n\!-\!2\delta\!-\!3s\!+\!6k\!-\!4)\!-\!\rho(\rho'\!-\!(n\!-\!2\delta
                   \!+\!2k \!-\!2))]\\
   &(\mbox{since $s\geq \delta+1$, $\rho'>n\!-\!(b\!+\!1)\delta\!-\!2$, $\delta\geq 2k$ and $k\geq 1$})\\
   &\geq\frac{(\delta\!-\!2k\!+\!1)(s\!-\!\delta)x_{1}y_{2}\rho'}{\rho(\rho'\!-\!(n\!-\!2\delta\!+\!2k\!-\!2))}\left(\frac{3n}{2}\!-\!4\delta\!+\!5k\!-\!3\!-\!\rho\right)~(\mbox{since $\rho\geq\rho'$ and $s\leq\frac{n+2k-2}{2}$})\\
    &=\frac{(\delta\!-\!2k\!+\!1)(s\!-\!\delta)x_{1}y_{2}\rho'}{\rho(\rho'\!-\!(n\!-\!2\delta\!+\!2k\!-\!2))}\left((n-1-\rho)+\frac{n}{2}+5k-4\delta-2\right)\\
    &>0~~(\mbox{since $\rho< n\!-\!1$, $n\geq 8\delta-10k+4$ and $s\geq \delta\!+\!1$}),
\end{aligned}
\end{equation*}
a contradiction. Therefore, we have $\rho'>\rho$, and it follows from  (\ref{equ::1}) and (\ref{equ::2}) that
$$\rho(G)\leq \rho(G^1)\leq \rho(G^2)<\rho(K_{\delta} \nabla (K_{n-2\delta+2k-1}\cup (\delta-2k+1)K_{1})).$$

{\flushleft\bf{Case 2.}} $s<\delta$.
\\[2mm]
\indent Let $G^{3}=K_{s} \nabla (K_{n-s-(\delta+1-s)(s-2k+1)}\cup (s-2k+1)K_{\delta+1-s})$. Recall that $G$ is a spanning subgraph of $G^1=K_{s} \nabla (K_{n_1}\cup K_{n_2}\cup \cdots \cup K_{n_{s-2k+2}})$, where $n_1\geq n_2\geq \cdots\geq n_{s-2k+2}$ and  $\sum_{i=1}^{s-2k+2}n_i=n-s$. Clearly, $n_{s-2k+2}\geq \delta+1-s$ because the minimum degree of $G^1$ is at least $\delta$. By Lemma \ref{lem::2.2}, we have
 \begin{equation}\label{equ::7}
\rho(G^{1})\leq \rho(G^{3}),
\end{equation}
where the equality holds if and only if $(n_1,\ldots,n_t)=(n-s-(\delta+1-s)(s-2k+1),\delta+1-s,\ldots,\delta+1-s)$.

If $s=2k$, then $G^3=K_{2k}\nabla(K_{n-\delta-1}\cup K_{\delta-2k+1})$. Combining Lemma \ref{lem::2.3}, (\ref{equ::1}) and (\ref{equ::7}), we get $$\rho(G)\leq \rho(G^1)\leq\rho(G^3)< \rho(K_{\delta}\nabla(K_{n-2\delta+2k-1}\cup (\delta-2k+1)K_{1})).$$
Thus, we consider $s\geq 2k+1$ in the following. Assume that $\rho(G^{3})=\rho^*\geq n -(\delta+1-s)(s-2k+1)$.
Let $x$ be the Perron vector of $A(G^{3})$. By symmetry, $x$ takes the same values  $x_{1}$, $x_{2}$, and $x_{3}$ on the vertices of $K_s$, $K_{\delta+1-s}$ and $K_{n-s-(\delta+1-s)(s-2k+1)}$, respectively. Then from $A(G^{3})x=\rho^* x$ we get
\begin{eqnarray}
\rho^* x_{1}&\!\!=\!\!&(s\!-\!1)x_1\!+\!(\delta\!+\!1\!-\!s)(s\!-\!2k\!+\!1)x_{2}\!+\!(n\!-\!s\!-\!(\delta\!+\!1\!-\!s)(s\!-\!2k\!+\!1))x_{3},\label{equ::8}\\
\rho^*x_{2}&\!\!=\!\!&sx_{1}+(\delta-s)x_{2}, \label{equ::9}\\
\rho^*x_{3}&\!\!=\!\!&sx_1+(n-s-1-(\delta+1-s)(s-2k+1))x_{3}.\label{equ::10}
\end{eqnarray}
Combining (\ref{equ::9}) and (\ref{equ::10}) yields that
\begin{equation}\label{equ::11}
\left\{
\begin{aligned}
 &x_{2}=\frac{sx_{1}}{\rho^*-\delta+s}, \\
 &x_{3}=\frac{sx_{1}}{\rho^*-(n-s-1-(\delta+1-s)(s-2k+1))}.
   \end{aligned}
   \right.
\end{equation}
Note that $n\geq F(k,\delta)\geq \delta(\delta-2k)^2+\delta-1$, $s\geq 2k+1$ and $\delta\geq s+1$. Then $\rho^*\geq n-(\delta+1-s)(s-2k+1)>\delta+1$. Combining this with (\ref{equ::8}) and (\ref{equ::11}), we have
\begin{equation*}
\begin{aligned}
&~~~~ \rho^*+1\\
&=s+\frac{s(\delta+1-s)(s-2k+1)}{\rho^*-\delta+s}+
\frac{s(n-s-(\delta+1-s)(s-2k+1))}{\rho^*-(n-s-1-(\delta+1-s)(s-2k+1))}\\
&< s+\frac{s(n-s)}{s+1} ~~( \mbox{since ~$\rho^*\geq n-(\delta+1-s)(s-2k+1)>\delta+1$})\\
 &=n\!-\!(\delta\!+\!1\!-\!s)(s-2k+1)\!-\!\frac{n\!-\!s\!-\!(\delta\!+\!1\!-\!s)(s-2k+1)(s\!+\!1)}{s\!+\!1}\\
 &\leq n\!-\!(\delta\!+\!1\!-\!s)(s-2k+1)\!-\!\frac{\delta(\delta-2k)^2+\delta-1\!-\!s\!-\!(\delta\!+\!1\!-\!s)(s-2k+1)(s\!+\!1)}{s\!+\!1}\\ &(\mbox{since ~$n\geq\delta(\delta-2k)^2+\delta-1$})\\
    &\leq n\!-\!(\delta\!+\!1\!-\!s)(s-2k+1)~~(\mbox{since ~$\delta\geq s+1$ and $s\geq 2k+1$})\\
     &\leq\rho^*,
 \end{aligned}
\end{equation*}
which is impossible. Therefore,
\begin{equation}\label{equ::12}
\begin{aligned}
                      \rho^*&<n-(\delta+1-s)(s-2k+1)\\
                         &=n-\delta+2k-2-((s-2k)(\delta-s)-1)\\
                         &\leq n-\delta+2k-2  ~~(\mbox{since ~$\delta\geq s+1$ and $s\geq2k+1$}).
\end{aligned}
\end{equation}
Since $K_{\delta}\nabla(K_{n-2\delta+2k-1}\cup(\delta-2k+1)K_1)$ contains $K_{n-\delta+2k-1}\cup(\delta-2k+1)K_1$ as a proper subgraph, we have $\rho(K_{\delta}\nabla(K_{n-2\delta+2k-1}\cup(\delta-2k+1)K_1))>n-\delta+2k-2$. Combining this with (\ref{equ::1}), (\ref{equ::7}) and (\ref{equ::12}), we may conclude that
$$\rho(G)\leq \rho(G^1)\leq \rho(G^3)<\rho(K_{\delta}\nabla(K_{n-2\delta+2k-1}\cup(\delta-2k+1)K_1)).$$

{\flushleft\bf{Case 3.}} $s=\delta.$
\\[2mm]
\indent By Lemma \ref{lem::2.2}, we have
$$\rho(G^1)\leq \rho(K_{\delta}\nabla(K_{n-2\delta+2k-1}\cup(\delta-2k+1)K_1)),$$
with equality holding if and only if $G^1\cong K_{\delta}\nabla(K_{n-2\delta+2k-1}\cup(\delta-2k+1)K_1)$. Combining this with (\ref{equ::1}), we conclude that
$$\rho(G)\leq\rho(K_{\delta}\nabla(K_{n-2\delta+2k-1}\cup(\delta-2k+1)K_1)),$$
where the equality holds if and only if $G\cong K_{\delta}\nabla(K_{n-2\delta+2k-1}\cup(\delta-2k+1)K_1)$. Observe that $K_{\delta}\nabla(K_{n-2\delta+2k-1}\cup(\delta-2k+1)K_1)$ is not $k$-extendable. Thus the result follows.\end{proof}

%From Theorem \ref{thm::1.1}, we can get the result of Liu, Liu and Feng\cite{W.L} immediately.
%\begin{cor}(See \cite{W.L})
%Let $G$ be a connected graph of even order $n$ with minimum degree $\delta(G)\geq 2$. If $n\geq \max\{7+7\delta+2\delta^{2}, \delta^{3}+3\delta^{2}+2\delta\}$ and $\rho(G)\geq\rho(K_{\delta}\nabla(K_{n-2\delta-1}\cup K_{\delta+1}))$,
%then $G$ has a perfect matching, unless $G\cong K_{\delta}\nabla(K_{n-2\delta-1}\cup K_{\delta+1})$.
%\end{cor}
%

By using Lemma \ref{lem::2.6}, we can give a short proof of Theorem \ref{thm::1.2}.

\renewcommand\proofname{\bf Proof of Theorem \ref{thm::1.2}}
\begin{proof}

Let $G=(A,B)$ be a balanced bipartite graph of order $n$ with minimum degree $\delta$. Suppose that $G$ is not $k$-extendable, by Lemma \ref{lem::2.4}, there exists some  nonempty subset $S\subseteq A$ with $s=|S|\leq |A|-k$ such that $|N(S)|<s+k$. Then $G$ is a spanning subgraph of $K_{s,s+k+1}\nabla_{1} K_{\frac{n}{2}-s,\frac{n}{2}-s-k-1}$ for some $s$ with $\delta \leq s\leq (n-2k-2)/4$. By Lemma \ref{lem::2.6}, we have
\begin{equation*}
\begin{aligned}
\rho(G)&\leq\rho(K_{s,s+k+1}\nabla_{1} K_{\frac{n}{2}-s,\frac{n}{2}-s-k-1})
\leq\rho(K_{\delta,\delta+k+1}\nabla_{1} K_{\frac{n}{2}-\delta,\frac{n}{2}-\delta-k-1}),
\end{aligned}
\end{equation*}
 where the first equality holds if and only if $G\cong K_{s,s+k+1}\nabla_{1} K_{\frac{n}{2}-s,\frac{n}{2}-s-k-1}$, and the second equality holds if and only if $s=\delta$. Note that $K_{s,s+k+1}\nabla_{1} K_{\frac{n}{2}-s,\frac{n}{2}-s-k-1}$ is not $k$-extendable. Thus the result follows.
\end{proof}

\section{Proof of Theorem~1.3}
\indent In this section, we give the proof of Theorem \ref{thm::1.3}. Before proceeding, the following structural lemma is needed. For any $v\in V(G)$, let $N_{G}(v)$ and $d_G(v)$ be the neighborhood and degree of $v$, respectively. For any subset $S\subset V(G)$, let $d_{S}(v)=|N_{G}(v)\cap S|$.
Let $f$ be a nonnegative integer-valued function on $V(G)$. An $f$-factor of $G$ is a spanning subgraph $H$ of $G$  such that $d_{H}(v)=f(v)$ for any $v\in V(G)$. It is clear that a $k$-factor is an $f$-factor with $f(v)=k$ for any $v\in V(G)$.
\begin{lem}[Ore \cite{Ore}]\label{lem::3.1}
A bipartite graph $G=(A,B)$ contains a $f$-factor if and only if $\sum_{x\in A}f(x)=\sum_{y\in B}f(y)$ and for any
subset $X\subseteq A$,
$$\sum_{x\in X}f(x)\leq \sum_{y\in N_{G}(X)}\min\{f(y), d_{X}(y)\}.$$
\end{lem}

\begin{lem}[Favaron, Mah\'{e}o and Sacl\'{e} \cite{O.F}]\label{lem::3.2}
If $G$ is a connected graph, then
$$\rho(G)\leq\max_{v\in V(G)}\sqrt{\sum_{u\in N_{G}(v)}d_{G}(u)},$$
with equality if and only if $G$ is either a regular graph or a semiregular bipartite graph.
\end{lem}

For any $S\subseteq V(G)$, let $N_{G}(S)$ be the set of neighbors of $v\in S$ in $G$. For $X,Y\subset V(G)$, we denote by $e(X)$ the number of edges in $G[X]$, and $e(X,Y)$ the number of edges with one endpoint in $X$ and one endpoint in $Y$. The following fact was pointed out in \cite{M.A}.
\begin{equation}\label{equ::13}
\begin{aligned}
\sum_{v\in N_{G}(u)}d_{G}(v)=d_{G}(u)+2e(N_{G}(u))+e(N_{G}(u),V(G)\backslash (N_{G}(u)\cup\{u\})).
\end{aligned}
\end{equation}

 Now we are in a position to give the proof of Theorems \ref{thm::1.3}.

\renewcommand\proofname{\bf Proof of Theorem \ref{thm::1.3}}
\begin{proof}
Let $G=(A,B)$ be a connected balanced bipartite graph of order $n$ with
$\rho(G)\geq \rho(K_{\frac{n}{2},\frac{n}{2}}\backslash E(K_{1,\frac{n}{2}-k+1}))$ and
$G\ncong K_{\frac{n}{2},\frac{n}{2}}\backslash E(K_{1,\frac{n}{2}-k+1})$. Suppose to the contrary that $G$ contains no $k$-factors, by Lemma \ref{lem::3.1}, there exists some nonempty subset $X\subseteq A$ with minimum cardinality such that
\begin{equation}\label{equ::14}
\begin{aligned}
k|X|> \sum_{y\in N_{G}(X)}\min\{k, d_{X}(y)\}.
\end{aligned}
\end{equation}
Let $Y=N_{G}(X)$. Then we partition $Y$ into $Y_{1}\cup Y_2$, where $Y_{1}=\{y\in Y|~d_{X}(y)\geq k\}$ and $Y_{2}=\{y\in Y|~d_{X}(y)<k\}$.
Suppose $|X|=x$ and $|Y_{i}|=y_{i}$ for $i=1,2$. Then, from (\ref{equ::14}), we can deduce that
\begin{equation*}
\begin{aligned}
kx> \sum_{y\in Y}\min\{k, d_{X}(y)\}=ky_{1}+\sum_{y\in Y_2}d_{X}(y).
\end{aligned}
\end{equation*}
Thus
\begin{equation}\label{equ::15}
\begin{aligned}
\sum_{y\in Y_2}d_{X}(y)<k(x-y_1),
\end{aligned}
\end{equation}
and
\begin{equation}\label{equ::16}
\begin{aligned}
x>y_{1}.
\end{aligned}
\end{equation}
We have the following three claims.

{\flushleft\bf Claim 1.} $d_{G}(v)\geq k$ for any $v\in V(G)$.\\[2mm]
\indent If not, without loss of generality, we assume that there exists a vertex $v\in A$ such that $d_{G}(v)<k$. Clearly, $G$ is a spanning subgraph of $K_{\frac{n}{2},\frac{n}{2}}\backslash E(K_{1,\frac{n}{2}-k+1})$. Then we have $\rho(K_{\frac{n}{2},\frac{n}{2}}\backslash E(K_{1,\frac{n}{2}-k+1}))\geq \rho(G)$, where the equality holds if and only if $G\cong K_{\frac{n}{2},\frac{n}{2}}\backslash E(K_{1,\frac{n}{2}-k+1})$. This is impossible because $\rho(G)\geq \rho(K_{\frac{n}{2},\frac{n}{2}}\backslash E(K_{1,\frac{n}{2}-k+1}))$ and $G\ncong K_{\frac{n}{2},\frac{n}{2}}\backslash E(K_{1,\frac{n}{2}-k+1})$. We complete the  proof of Claim 1.\qed

Take $u\in V(G)$ with $R_{u}=\max_{v\in V(G)}R_{v}$.

{\flushleft\bf Claim 2.} $R_{u}>\frac{n}{2}\left(\frac{n}{2}-1\right)$.
\\[2mm]
\indent Since $K_{\frac{n}{2}-1,\frac{n}{2}}$ is a proper subgraph of $K_{\frac{n}{2},\frac{n}{2}}\backslash E(K_{1,\frac{n}{2}-k+1})$ for $k\geq 2$, we have
$$\rho(G)\geq \rho(K_{\frac{n}{2},\frac{n}{2}}\backslash E(K_{1,\frac{n}{2}-k+1}))>\sqrt{\frac{n}{2}\left(\frac{n}{2}-1\right)}.$$
By Lemma \ref{lem::3.2}, we can deduce that
$$R_{u}\geq \rho^{2}(G)\geq  \rho^{2}(K_{\frac{n}{2},\frac{n}{2}}\backslash E(K_{1,\frac{n}{2}-k+1}))>\frac{n}{2}\left(\frac{n}{2}-1\right).$$
This completes the proof of Claim 2.\qed

{\flushleft\bf Claim 3.} $x\geq 2$.
\\[2mm]
\indent If not, suppose that $X=\{w\}$, then $d_{X}(y)=1$ for $y\in N_{G}(w)$, and hence
$$\sum_{y\in N_{G}(w)}\min\{k, d_{X}(y)\}=d_{G}(w)<k$$
by (\ref{equ::14}) and $k\geq 2$, which contradicts Claim 1. Therefore, $x\geq 2$. This completes the proof of Claim 3.\qed

Now we divide the proof into two situations.
{\flushleft\bf Case 1.} $u\in A$.
\\[2mm]
\indent Note that $N_{B}(u)\subseteq B$ and $e(N_{B}(u))=0$. By (\ref{equ::13}), we have
\begin{equation}\label{equ::17}
\begin{aligned}
R_{u}=d_{B}(u)+e(N_{B}(u),A\backslash\{u\}).
\end{aligned}
\end{equation}
We assert that $d_{B}(u)=|B|=n/2$. If not, $d_{B}(u)\leq n/2-1$.
Observe that $e(N_{B}(u),A\backslash\{u\})\leq \left(n/2-1\right)^{2}$.
Then, from Claim 2 and (\ref{equ::17}), we get
\begin{equation*}
\begin{aligned}
\frac{n}{2}\left(\frac{n}{2}-1\right) <R_{u}=d_{B}(u)+e(N_{B}(u),A\backslash\{u\})\leq \frac{n}{2}-1+\left(\frac{n}{2}-1\right)^{2}=\frac{n}{2}\left(\frac{n}{2}-1\right),
\end{aligned}
\end{equation*}
a contradiction. Thus $N_{B}(u)=B$.

If $u\notin X$, then $2\leq x\leq n/2-1$ by Claim 3. We consider the following two subcases.

{\flushleft\bf Subcase 1.1.} $2\leq x\leq k$.
\\[2mm]
\indent Since $2\leq x\leq k$, we have $d_{X}(y)\leq x\leq k$ for $y\in Y$,
and hence
$$kx>\sum_{y\in Y}\min\{k, d_{X}(y)\}=\sum_{y\in Y}d_{X}(y)$$
by (\ref{equ::14}). For $w\in X$, let $X'=X-w$ and $Y'=N_{B}(X')$. Then $d_{X'}(y)\leq d_{X}(y)\leq k$ for $y\in Y$. Combining this with the minimality of $X$, we obtain
$$k(x-1)=k|X'|\leq \sum_{y\in Y'}\min\{k, d_{X'}(y)\}= \sum_{y\in Y'}d_{X'}(y)=\sum_{y\in Y}d_{X}(y)-d_{Y}(w)< kx-d_{Y}(w),$$
and hence $d_{Y}(w)=d_{G}(w)<k$, which contradicts Claim 1.

{\flushleft\bf Subcase 1.2.} $k+1\leq x\leq n/2-1$.
\\[2mm]
\indent Recall that $X\subseteq A$ and $Y=N_{G}(X)$. Then $e(X, B-Y)=0$. Combining this with Claim 2, (\ref{equ::15}) and (\ref{equ::16}), we have
\begin{equation*}
\begin{aligned}
\frac{n}{2}\left(\frac{n}{2}-1\right)&<R_{u}\\
&=d_{B}(u)+e(B,A\backslash\{u\})\\
&=e(A,B)\\
&=e(X,Y_1)+e(X,Y_2)+e(A-X,B)\\
&<xy_1+k(x-y_1)+\frac{n}{2}\left(\frac{n}{2}-x\right) ~~(\mbox{by (\ref{equ::15}}))\\
&=\frac{n}{2}\left(\frac{n}{2}-1\right)+kx+(x-k)y_1-\frac{n}{2}(x-1)\\
&\leq \frac{n}{2}\left(\frac{n}{2}-1\right)+kx+(x-k)(x-1)-\frac{n}{2}(x-1) ~~(\mbox{by (\ref{equ::16}}))\\
&=\frac{n}{2}\left(\frac{n}{2}-1\right)+\left(x-\frac{n}{2}\right)(x-1)+k\\
&\leq \frac{n}{2}\left(\frac{n}{2}-1\right)+\left(x-\frac{n}{2}\right)(x-1)+x-1~~(\mbox{since $x\geq k+1$})\\
&= \frac{n}{2}\left(\frac{n}{2}-1\right)+\left(x-\frac{n}{2}+1\right)(x-1)\\
&\leq \frac{n}{2}\left(\frac{n}{2}-1\right)~~(\mbox{since $k+1\leq x\leq \frac{n}{2}-1$ and $k\geq 2$}),
\end{aligned}
\end{equation*}
a contradiction.

If $u\in X$, since $N_{Y}(u)=B$, we have $Y=B$. For $2\leq x \leq n/2-1$, we  also can deduce a contradiction by using a similar analysis as above. For $x=n/2$, we have $X=A$. We assert that $Y_2=\emptyset$. If not, there exists some vertex $z\in Y_2$ such that $d_{G}(z)=d_{A}(z)<k$, which contradicts Claim 1. This implies that $Y_1=B$, and hence $y_1=n/2$. Combining this with $d_{X}(y)\geq k$ for $y\in Y_1$ and (\ref{equ::14}), we have
$$\frac{nk}{2}=kx>\sum_{y\in Y}\min\{k,d_{X}(y)\}=ky_1=\frac{nk}{2},$$
a contradiction.

{\flushleft\bf Case 2.} $u\in B$.
\\[2mm]
\indent If $2\leq x\leq n/2-1$, by using a similar analysis as above, we can deduce a contradiction. If $x=n/2$, then $X=A$. First we assert that $B=Y$. If not, $B-Y\neq \emptyset$. Since $e(X,B-Y)=0$,  $G$ is disconnected. This is impossible because $G$ is a connected balanced bipartite graph. Next we assert that $Y_2=\emptyset$. If not, there exists some vertex $z\in Y_2$ such that $d_{G}(z)=d_{A}(z)<k$, which contradicts Claim 1. This implies that $Y_1=B$ and $y_1=n/ 2$. Combining this with $d_{X}(y)\geq k$ for $y\in Y_1$ and (\ref{equ::14}), we have
$$\frac{nk}{2}=kx>\sum_{y\in Y}\min\{k,d_{X}(y)\}=ky_1=\frac{nk}{2},$$
which is also impossible.
\end{proof}

\section{Concluding remarks}

\begin{lem}[Liu,  Qian, Sun and Xu \cite{L-G}]\label{lem::4.1}
If $G$ is a bipartite graph with $m\geq 1$ edges and $n$ vertices, then
$$\rho(G)\leq \sqrt{m},$$
and equality holds if and only if $G\cong K_{p,q}\cup (n-p-q)K_1$, where $pq=m$.
\end{lem}

In 2012, Lu, Liu and Tian\cite{Lu-Liu} investigated the condition of $e(G)$ for the existence of a hamiltonian cycle in a graph, where $e(G)$ denote the number of edges in $G$.

\begin{lem}[Lu, Liu and Tian \cite{Lu-Liu}]\label{lem::4.2}
Let $G$ be a balanced bipartite graph with order $n\geq 8$ and minimum degree $\delta\geq 1$. If
$$e(G)\geq \frac{n}{2}\left(\frac{n}{2}-1\right)+1,$$
then $G$ contains a hamiltonian cycle, unless $G\cong K_{\frac{n}{2},\frac{n}{2}}\backslash E(K_{1,\frac{n}{2}-1})$.
\end{lem}

Utilizing Lemmas \ref{lem::4.1} and \ref{lem::4.2}, we can easily deduce the following result.
\begin{thm}\label{lem::4.3}
Let $G$ be a balanced bipartite graph with order $n\geq 8$. If
$$\rho(G)\geq \rho(K_{\frac{n}{2},\frac{n}{2}}\backslash E(K_{1,\frac{n}{2}-1})),$$
then $G$ contains a hamiltonian cycle, unless $G\cong K_{\frac{n}{2},\frac{n}{2}}\backslash E(K_{1,\frac{n}{2}-1})$.
\end{thm}
\renewcommand\proofname{\bf Proof}
\begin{proof}
Note that $K_{\frac{n}{2},\frac{n}{2}-1}$ is a proper subgraph of $K_{\frac{n}{2},\frac{n}{2}}\backslash E(K_{1,\frac{n}{2}-1})$. Combining this with Lemma \ref{lem::4.1}, we have
$$ \sqrt{\frac{n}{2}\left(\frac{n}{2}-1\right)}<\rho(G)\leq\sqrt{e(G)},$$
which leads to $e(G)\geq \frac{n}{2}\left(\frac{n}{2}-1\right)+1$. Thus, the result follows by Lemma \ref{lem::4.2}.
\end{proof}
Note that a hamiltonian cycle is a connected $2$-factor. It is interesting to see that Theorem \ref{thm::1.3} also provides a spectral condition for the existence of a connected $2$-factor in a balanced bipartite graph. Naturally, we have the following problem.
\begin{prob}
For $k\geq 3$. Suppose that $G$ is a balanced bipartite graph of order $n$ with connected $k$-factors. Does $\rho(G)\geq \rho(K_{\frac{n}{2},\frac{n}{2}}\backslash E(K_{1,\frac{n}{2}-k+1})))$, unless $G\cong K_{\frac{n}{2},\frac{n}{2}}\backslash E(K_{1,\frac{n}{2}-k+1})$?
\end{prob}

A graph $G$ is said to be \textit{$k$-factor-critical}, if $G-S$ has a perfect matching for every subset $S\subseteq V(G)$ with $|S|=k$. The following fundamental lemma provides a sufficient and necessary condition for a graph to be $k$-factor-critical was obtained by Favaron\cite{Favaron} and Yu\cite{Yu}, independently.
\begin{lem}[Favaron\cite{Favaron}; Yu\cite{Yu}]\label{lem::4.4}
For $k\geq 1$, a graph $G$ of order $n$ is $k$-factor-critical if and if $n\equiv k~(\rm{mod}~2)$ and
$$o(G-S)\leq |S|-k$$
for any subset $S\subseteq V(G)$ with $|S|\geq k$.
\end{lem}

Different types of matching extensions are closely related. It is clear that a $2k$-factor-critical graph must be $k$-extendable. On the other hand, not all $k$-extendable graphs are $2k$-factor-critical. For example, $K_{k+1,k+1}$ is $k$-extendable but is not $2k$-factor-critical.
By Lemma \ref{lem::4.4} and by using a similar analysis as in the proof of Theorem \ref{thm::1.1}, we can easily deduce the following result.

\begin{thm}
Suppose that $G$ is a connected graph of order $n\geq\max\{8\delta-5k+4, \delta(\delta-k)^2+\delta-1\}$ with minimum degree $\delta\geq k$, where $n\equiv k~(\rm{mod}~2)$ and $k\geq 1$. If $$\rho(G)\geq\rho(K_{\delta}\nabla(K_{n-2\delta+k-1}\cup(\delta-k+1)K_1)),$$ then $G$ is $k$-factor-critical, unless $G\cong K_{\delta}\nabla(K_{n-2\delta+k-1}\cup(\delta-k+1)K_1)$.
\end{thm}

\end{document}